\documentclass[12pt]{amsart}
\usepackage{amsmath,amsthm,amsfonts,amssymb,times, epsfig}
\DeclareMathAlphabet{\curly}{U}{rsfs}{m}{n}
\theoremstyle{remark}

\theoremstyle{plain}

\newtheorem{lem}{Lemma}[section]
\newtheorem{thm}{Theorem}
\newtheorem{cor}{Corollary}
\numberwithin{equation}{section}

\newcommand{\Z}{{\mathbb Z}}

\newcommand{\be}{\begin{equation}}
\newcommand{\ee}{\end{equation}}
\newcommand{\benn}{\begin{equation*}}
\newcommand{\eenn}{\end{equation*}}

\renewcommand{\a}{\ensuremath{\alpha}}

\newcommand{\del}{\ensuremath{\delta}}
\newcommand{\eps}{\ensuremath{\varepsilon}}
\renewcommand{\(}{\left(}
\renewcommand{\)}{\right)}

\newcommand{\pfrac}[2]{\left(\frac{#1}{#2}\right)}
\newcommand{\sfrac}[2]{{#1/#2}}

\renewcommand{\th}{\theta}
\newcommand{\li}{{\rm li}}

\begin{document}
\baselineskip 17pt
\title{Divisors of the Euler and Carmichael functions}
\date{\today}
\author{Kevin Ford and Yong Hu}

\begin{abstract}
We study the distribution of divisors of Euler's totient function and
Carmichael's function. In particular, we estimate how often the values of these
functions have "dense" divisors.
\end{abstract}

\address{Department of Mathematics, 1409 West Green St.,
University of Illinois, Urbana, IL 61801, USA}

\email{KF: ford@math.uiuc.edu; YH: cofuns@gmail.com}

\thanks{First author supported in part by National Science Foundation
  grant DMS-0555367}

\maketitle

\section{Introduction}
Two of the most studied functions in the theory of numbers are Euler's
totient function $\phi(n)$ and Carmichael's function $\lambda(n)$, the
first giving the order of the group $(\Z/n\Z)^*$ of reduced residues
modulo $n$, and the latter giving the maximum order of any element of
$(\Z/n\Z)^*$. The distribution of $\phi(n)$ and $\lambda(n)$ has been
investigated from a variety of perspectives. In particular, many
interesting properties of these functions require knowledge of the
distribution of prime factors of $\phi(n)$ and $\lambda(n)$, e.g.,
\cite{Erdosphi}, \cite{EP}, \cite{EGPS}, \cite{EPS}, \cite{totients},
\cite{LP}, \cite{Warlimont}.

The distribution of all of the divisors of $\phi(n)$ and $\lambda(n)$
has thus far received little attention, perhaps due to the complicated
way in which prime factors interact to form divisors. From results
about the normal number of prime factors of $\phi(n)$ and $\lambda(n)$
\cite{EP}, one deduces immediately that $\tau(\phi(n))$ and
$\tau(\lambda(n))$ are each $\exp\{ \frac{\log 2}{2} (\log\log n)^2\}$
for almost all $n$. However, the determination of the \emph{average}
size of $\tau(\phi(n))$ and of $\tau(\lambda(n))$ is more complex, and
has been studied recently by Luca and Pomerance \cite{LP07}.

In this note we investigate problems about localization of divisors of
$\phi(n)$ and $\lambda(n)$. Our results have application to the
structure of $(\Z/n\Z)^*$, since the set of divisors of $\lambda(n)$
is precisely the set of orders of elements of $(\Z/n\Z)^*$. We say
that a positive integer $m$ has $u$-dense divisors (for short, $m$ is
\emph{$u$-dense}) 
if whenever $1\le y<m$, there is a divisor of $m$ in the interval
$(y,uy]$. The distribution of $u$-dense numbers for general $u$ has
been investigated by Tenenbaum (\cite{Ten86}, \cite{Ten95}) and Saias
(\cite{Saias1}, \cite{Saias2}). According to Th\'eor\`eme 1 of
\cite{Saias1}, the number of $u$-dense integers $m \le x$ is $\asymp
\sfrac{(x\log u)}{\log x}$, uniformly for $2 \le u \le x$. In
particular, the number of $2$-dense integers $m\le x$ is $\asymp
\sfrac{x}{\log x}$, that is, the $2$-dense integers are about as
sparse as the primes.

By contrast, we show that $2$-dense values of $\phi(n)$ and
$\lambda(n)$ are very common.

\begin{thm}\label{dense}
If $x$ is sufficiently large, then for $\gg x$ integers $n\le x$, both $\phi(n)$ and $\lambda(n)$ are $2$-dense.
\end{thm}

There are relatively simple heuristic reasons for believing Theorem
\ref{dense}. Recall that
\begin{align*}\label{philam}
\phi(p_1^{e_1} \cdots p_k^{e_k}) &= p_1^{e_1-1}(p_1-1) \cdots
p_k^{e_k-1}(p_k-1), \\ \lambda(p_1^{e_1} \cdots p_k^{e_k}) &=
\text{lcm}[\lambda(p_1^{e_1}),\ldots, \lambda(p_k^{e_k})],
\end{align*}
where $\lambda(p_i^{e_i})=\phi(p_i^{e_i})$ if $p_i$ is odd or $p_i=2$
and $e_i \le 2$, and $\lambda(2^e)=2^{e-2}$ for $e\ge 3$. In
particular, $\phi(n)$ and $\lambda(n)$ have the same prime factors. 
We recall that
\begin{itemize}
\item
most of these prime factors are factors of shifted primes $p-1$ where $p | n$
\item
for most primes $p$, $p-1$ has about $\log \log p$ prime factors 
\cite{Erdosphi}
\item
for most integers $n$, 
we have $\omega(n)$ about $\log\log n$ and, writing these 
distinct prime factors as $p_1(n) < p_2(n) < \dots < p_{\omega(n)}(n)$,
we have for all $k$ except for the smallest ones
$\log \log p_k(n) \approx k$ (see, e.g., Theorem 10 of \cite{divisors})
\item
$\log \phi(n) \sim \log n$.
\end{itemize}
With these four facts, we find that for most values of $n$
\begin{align*}
\Omega (\phi(n)) &\approx \sum_{p|n} \Omega(p-1) \approx
\sum_{p|n} \log \log p \\ 
&\approx \sum_{k\leq \log\log n} k \sim
(\log \log n)^2 / 2 \sim (\log \log \phi(n))^2 / 2.
\end{align*}
See \cite{EP} for a more precise result about the normal behavior of
$\Omega (\phi(n))$.
On the other hand, for most values of $m$, $\Omega(m)$ is about
$\log\log m$. 
So, usually $\phi(n)$ has far
more divisors than a typical integer of its size.


We therefore expect the divisors of $\phi(n)$, especially the smaller
divisors, to be ``very dense'' for most $n$, and the same should be
true of small divisors of $\lambda(n)$. On the other hand, there are a
large proportion of $n$ for which the divisors of $\phi(n)$ and
$\lambda(n)$ are not very dense. To state our next result, we define
$\theta$ to be the supremum of real numbers $c$ so that there are $\gg
x/\log x$ primes $p \le x$ with $p-1$ having a prime factor
$>p^c$. Many papers have been written on bounding $\th$, and the
current record is $\th \ge 0.677$ and due to Baker and Harman
\cite{BH}.

\begin{thm}\label{nondense}
Let $0 <c < 2 \th-1$. If $x$ is sufficiently large, then for $\gg_c x$ of the integers $n\le x$, neither $\phi(n)$ nor $\lambda(n)$ is $x^c$-dense.
\end{thm}

It is conjectured that $\th=1$, and this would imply the conclusion of Theorem \ref{nondense} for any $c<1$.

If $u < 2$, there are no $u$-dense integers $m>1$. However, it is
possible that the divisors of a given integer in some long interval do
have consecutive ratios which are $\le u$. We say that an integer $n$
has $u$-dense divisors in a set $I$ (for short, $n$ is $u$-dense in a
set $I$) if for every $y \in I$, the interval
$(y,uy]$ contains a divisor of $n$. The following makes precise what
we claimed earlier about the ``very dense'' nature of the small
divisors of $\phi(n)$ and $\lambda(n)$.

\begin{thm}\label{verydense}
For every positive integer $h$ and $0 < \delta < 1$, there is a constant $c = c(h, \del) > 0$ so that if $x$ is sufficiently large, then for more than $(1 - \delta)x$ of the integers $n \le x$, $\phi(n)$ and $\lambda(n)$ are both $(1+1/h)$-dense in $[h, x^c]$.
\end{thm}

Notice that the left endpoint $h$ of the interval cannot be replaced by $h-1$, since if $h - 1 \le a < \sfrac{h}{(1 + \sfrac{1}{h})}$, there are no integers in $(a, a(1 + \sfrac{1}{h})]$. Likewise, if we assume that $\th=1$, then we cannot take $c$ independent of $\delta$ in light of Theorem \ref{nondense}.

Using Theorem \ref{verydense}, we prove a more general version of Theorem \ref{dense}.

\begin{thm}\label{denseh}
For every positive integer $h$, there are $\gg_h x$ integers $n\le x$ such that $\phi(n)$ is $(1+1/h)$-dense in $[h, \phi(n)/(h + 1))$ and $\lambda(n)$ is $(1+1/h)$-dense in $[h, \lambda(n)/(h + 1))$.
\end{thm}

We also record a limiting case of Theorem \ref{verydense}.

\begin{cor}\label{limitverydense}
Suppose $g(x)$ is a positive function decreasing monotonically to 0 and let $h$ be a positive integer. Almost all $n\le x$ have the property that $\phi(n)$ and $\lambda(n)$ are $(1+1/h)$-dense in $[h,x^{g(x)}]$.
\end{cor}

Analogous to the problems studied in \cite{Hxyz}, \cite{Hxy2y}, \cite{Ten84}, we can study the distribution of integers with $\phi(n)$ having a divisor in a \emph{single} interval. Let
$$
B(x,y,z) = |\{n \le x: \exists \, d|\phi(n), y<d\le z \}|.
$$
An almost immediate corollary of Theorems \ref{dense}, \ref{nondense} and \ref{verydense} is the following result in the special case $z=2y$.

\begin{cor}\label{bxy2y}
\mbox{}
\begin{itemize}
\item[(i)] Uniformly for $1\le y\le x/2$, we have $B(x,y,2y) \gg x$.
\item[(ii)] Fix $1-\theta<c<1/2$. Then, uniformly for $x^c \le y\le x^{1-c}$,
we have $x-B(x,y,2y) \gg x$.
\item[(iii)] Let $g(x)\to 0$ monotonically. Then, for $1\le y\le x^{g(x)}$,
we have $B(x,y,2y)\sim x$.
\end{itemize}
\end{cor}

We leave as an open problem the determination of the order of magnitude of $B(x,y,z)$ for all $x,y,z$.

We note that easy modifications of our proofs give the same results for the sum of divisors function $\sigma(n)$ in place of $\phi(n)$, since $\sigma(p)=p+1$ for primes $p$.

The authors would like to thank Igor Shparlinski for posing the question to study the divisors of $\phi(n)$.

%
%
\section{Preliminaries}
%
%

Throughout this paper, the letters $p$ and $q$, with or without subscripts,
will always denote primes. Constants implied by the $O$ and $\ll$ symbols
are absolute, unless dependence on a parameter is indicated 
by a subscript. All constants are effectively computable as well.
We denote by $P^+(m)$ the largest prime factor of $m$, with the convention
that $P^+(1)=0$.

Our key lemma, presented below, says roughly
that the small \emph{prime} factors of $\phi(n)$ are quite dense.

\begin{lem}\label{Phinodivg}
For some large constant $C$, if $\frac{C}{\log x} \le g \le \frac{1}{10}$
and $\frac{1}{(g\log x)} \le \eps \le \frac14$, then the number of
$n\le x$ for which $\phi(n)$ does not have a prime divisor in
$(x^g,x^{g(1+\eps)}]$ is $\ll g^{\eps/2} \log(1/g) x$.
\end{lem}

\begin{proof}
First, we claim that for large $x$ and $w\ge x^{6g}$, that 
\be\label{divp-1} 
| \{ p \le w : p - 1 \mbox{ has no prime factor in } (x^g,x^{g(1+\eps)}]\}
| \le \(1 - \frac{2\eps}{3}\)\frac{w}{\log w}. 
\ee 
Let $\pi(w;r,a)$ be the number of primes $p\le w$ which satisfy $p\equiv a\pmod{r}$. For positive integer $r$, write
$$
\pi(w;r,1) = \frac{\li(w)}{\phi(r)} + E(w;r),
$$
where
$$
\li(w) = \int_2^w \frac{dt}{\log t}.
$$
Using the Bombieri-Vinogradov Theorem (\cite{Dav}, Ch. 28) and the
Mertens' estimates, the number of primes $p \le w$ such that $p - 1$
has a prime factor in $(x^g,x^{g(1+\eps)}]$ is
\begin{align*}
&\ge \sum_{x^g < q \le x^{(1 + \eps)g}}\pi(w;q,1) - \sum_{x^g
 < q_1 < q_2 \le x^{(1 + \eps)g}}\pi(w;q_1q_2,1) \\
&= \sum_q \( \frac{\li(w)}{q - 1} + E(w;q) \) -
 \sum_{q_1, q_2} \( \frac{\li(w)}{(q_1 - 1)(q_2 - 1)} + E(w;q_1q_2) \) \\
&= \li(w) \left[\log(1 + \eps) - \frac12 \log^2 (1 + \eps) + O\pfrac{1}{\log^2
 x^g}\right] + O\pfrac{w}{\log^3 w}\\ 
&\ge \frac{3\eps}{4}\frac{w}{\log w}.
\end{align*}
For the last step, we used the fact that $w \ge x^{6g} \ge e^{6C}$ and $C$ is sufficiently large. This proves \eqref{divp-1}.

Consider $x/\log x < n \le x$ such that $\phi(n)$ does not have a prime divisor in $(x^g,x^{g(1+\eps)}]$. We can write $n = q^{\a_1}_1q^{\a_2}_2...q^{\a_k}_km$, where $q_1 > q_2 > \dots > q_k > x^{6g}$, $\a_i \ge 1$ for $1 \le i \le k$ and $P^{+}(m) \le x^{6g}$. Then $q_1,\ldots,q_k\in T$, the set of primes $p$ such that $p - 1$ does not have a prime factor in $(x^g,x^{g(1+\eps)}]$. By \eqref{divp-1} and partial summation,
\begin{align*}
\sum_{\substack {x^{6g} < q \le x \\ q \in T}} \sum_{a\ge 1} \frac1{q^a} 
&\le \(1 -\frac{2\eps}{3}\)\( \log \frac{1}{6g} + \frac{1}{\log x} \) 
+ \sum_{q>e^{6C}} \frac{1}{q(q-1)} \\
&\le \(1 - \frac{\eps}{2} \) \log \frac{1}{6g}
\end{align*}
for sufficiently large $x$. By Theorem 07 of \cite{divisors}, for some positive constant $c_0$ and uniformly in $x\ge z$, $y\ge 2$, the number of integers $n\le x$ divisible by a number $m>z$ with $P^+(m)\le y$ is $\ll x \exp\{ -c_0 \frac{\log z}{\log y} \}$. Consequently, the number of $n$ with $m>x^{1/3}$ is $\ll x e^{-c_0/18g} \ll gx$. For other $n$, we may assume that $m \le x^{1/3}$, and thus $k \ge 1$. Again by the above theorem, the number of $n$ with $q_1 \le \log^{10} x$ is $\ll x/\log x \ll gx$. For remaining $n$, we have $q_1^{\alpha_1-1} \cdots q_k^{\alpha_k-1} \le \log^2 x$, for otherwise, $q_1^{\lfloor \sfrac{\alpha_1}{2} \rfloor} \cdots q_k^{\lfloor \sfrac{\alpha_k}{2} \rfloor} \ge q_1^{\sfrac{(\alpha_1-1)}{2}} \cdots q_k^{\sfrac{(\alpha_k-1)}{2}} > \log x$ and the number of $n$ divisible by $d^2$ for some $d> \log x$ is $O(x/\log x)$. Hence $q_1 \cdots q_k \ge x^{1/2}$. In particular, $q_1 \ge \max(x^{\frac{1}{2k}},\log^{10} x)$ and $\a_1=1$. Given $q_2^{\a_2}, \ldots, q_k^{\a_k}$, and $m$, the number of $q_1$ is, by the Chebyshev estimates for primes,
$$
\ll \frac{x}{q_2^{\a_2} \cdots q_k^{\a_k} m\log(x/(q_2^{\a_2}\cdots q_k^{\a_k}
 m))} \ll \frac{kx}{\log x} \, \frac{1}{q_2^{\a_2} \cdots q_k^{\a_k} m}.
$$
Given $q_2^{\a_2},\ldots,q_k^{\a_k}$,
$$
\sum_{P^{+}(m) \le x^{6g}} \frac1m \ll \log (x^{6g}) = 6g \log x.
$$
With fixed $k$, we have
\begin{align*}
\sum_{q_2, \ldots, q_k\in T} \sum_{\a_2,\ldots,\a_k\ge 1} 
\frac{1}{q_2^{\a_2} \cdots q_k^{\a_k}} &\le
\frac{1}{(k - 1)!}\biggl(\sum_{\substack {x^{6g} < q \le x \\ q \in T}}
\sum_{a\ge 1}\frac1{q^a} \biggr)^{k - 1} \\
&\le \frac{((1- \frac{\eps}{2})\log \frac{1}{6g})^{k - 1}}{(k - 1)!}.
\end{align*}
The total number of such $n$ is
\begin{align*}
&\ll gx + gx \sum_{1\le k\le 1/(6g)} \frac{k}{(k - 1)!}\(\(1-
\frac{\eps}{2}\)\log \frac{1}{6g}\)^{k - 1} \\
&\ll gx + gx \( \log \frac{1}{6g} \) \sum_{j = 0}^{\infty} \frac{((1-
\frac{\eps}{2})\log \frac{1}{6g})^{j}}{j!} \\
&= gx + gx\(\log \frac{1}{6g}\)\(\frac{1}{6g}\)^{1-
\frac{\eps}{2}} \\
&\ll g^{\frac{\eps}{2}} \(\log \frac{1}{g}\)x.
\end{align*}
This completes the proof.
\end{proof}

{\bf Remarks.} Since $\phi(n)$ and $\lambda(n)$ have the same prime
 factors, Lemma \ref{Phinodivg} holds with $\phi$ replaced by 
$\lambda$. With a finer 
analysis, it is possible to remove the factor $\log(1/g)$ appearing
 in the conclusion of Lemma \ref{Phinodivg}. Also, if $\eps$ is fixed,
 then $g^{\eps/2}\log (1/g) \ll_\eps g^{\eps/3}$, an inequality we
 shall
 use in the application of Lemma \ref{Phinodivg}.

We next give a method of constructing integers which are dense in an interval.

\begin{lem}\label{densecondinterval}
Suppose that $h$ is a positive integer, $y\ge h$, and $D$ is
 $(1+1/h)$-dense in $[h,y]$. Suppose also that $m=Dm_1\cdots m_k$,
 where for $1\le j\le k$, $m_j \le (y/h) m_1\cdots m_{j-1}$. Then 
$m$ is $(1+1/h)$-dense in $[h,m_1\cdots m_k y]$.
\end{lem}

\begin{proof}
By hypothesis, the lemma holds for $k=0$. Suppose the lemma is true
 for $k=l$, $m$ satisfies the hypotheses with $k=l+1$ and put
 $m'=Dm_1\cdots m_{l}$. Then $m'$ is $(1+1/h)$-dense in $[h,m_1\cdots
 m_l y]$. Multiplying the divisors of $m'$ by $m_{l+1}$, we find that
 $m$ is also $(1+1/h)$-dense in $[m_{l+1} h,m_1\cdots m_{l+1}y]$. Our
 assumption about $m_{l+1}$ implies that $m$ is $(1+1/h)$-dense in 
$[h,m_1\cdots m_{l+1}y]$, as desired.
\end{proof}

\begin{lem}\label{smalldiv}
Given any positive integer $D$, $n$ is divisible by a prime $q \equiv 1\pmod{D}$ for almost all $n$.
\end{lem}

\begin{proof}
By a theorem of Landau \cite{La}, the number of $n\le x$ which have no prime factor $q\equiv 1\pmod{D}$ is asymptotic to $c(D) x (\log x)^{-1/\phi(D)}$ for some constant $c=c(D)$.
\end{proof}

Luca and Pomerance \cite{LP} have recently proven a stronger
statement, namely that for some constant $c_1$, for almost all 
integers $n$, $\phi(n)$ is divisible by every prime power 
$\le c_1 \frac{\log\log n}{\log\log\log n}$.

%
%
\section{Proof of the theorems}
%
%

\begin{proof}[Proof of Theorem \ref{verydense}]
Fix $h$ and $\delta$, and let $y$ be sufficiently large, depending on 
$h$, and such that $y > h^5$. Let $D$ be the product of all prime 
powers $\le y$. Let $\eps=\frac14$ and let $Y=(y/h)^{4/5}$. Let $C$ 
be the constant in Lemma \ref{Phinodivg}.

Consider the intervals $I_j=(Y^{(5/4)^{j-1}}, Y^{{(5/4)}^{j}}]$ 
$(1\le j\le J)$, where $Y\ge e^C$. Fix $c$ so that $0 < c \le 
\sfrac{1}{20}$, let $x$ be sufficiently large, and take $J$ so 
that $Y^{(5/4)^{J-2}} < x^c \le Y^{(5/4)^{J-1}}$. Then $Y^{(5/4)^J}
 < (Y^{(5/4)^{J-2}})^2 < x^{2c} \le x^{1/10}$. By Lemma \ref{Phinodivg}, 
if $y$ is large enough, then the number of integers $n \le x$ for 
which $\phi(n)$ does not have prime factors in $I_j$ is
$$
\ll \(\frac{\log Y^{(5/4)^{j-1}}}{\log x}\)^{\sfrac{1}{12}} x.
$$
Summing over $j$, we find that $\phi(n)$ has a prime factor in 
every interval $I_j$ for all $n\le x$ except for a set of size
$$
\ll \pfrac{\log Y^{(5/4)^J}}{\log x}^{\sfrac{1}{12}} x < (2c)^{\sfrac{1}{12}}x.
$$

If $c$ is small enough, for at least $(1-\delta/2)x$ of the integers $n\le x$, $\phi(n)$ has a prime factor in every interval $I_j$. Applying Lemma \ref{smalldiv}, for at least $(1- \delta) x$ integers $n\le x$, $\phi(n)$ is divisible by a prime $q \equiv 1\pmod{D^3}$ and has a prime factor in every interval $I_j$. For each such $n$, let $p_1,\ldots,p_J$ be primes dividing $\phi(n)$ and such that $p_j\in I_j$ for $1\le j\le J$. By hypothesis, $p_3>(y/h)^{5/4}>y$, hence $p_j\nmid D$ for $j\ge 3$. Since $D^3|(q-1)|\lambda(n)|\phi(n)$, we have that $\lambda(n)$ and $\phi(n)$ are each divisible by $Dp_1\cdots p_J$. By definition, $D$ is divisible by every positive integer $\le y$, hence $D$ is $(1+1/h)$-dense in $[h,y]$. Also, $p_1 \le Y^{5/4} = y/h$, and for $j\ge 2$,
$$
p_j \le Y^{(5/4)^j} \le Y^{5/4} \prod_{1\le i\le j-1} Y^{(5/4)^{i-1}}
\le (y/h) p_1\cdots p_{j-1}.
$$
By Lemma \ref{densecondinterval}, $\phi(n)$ and $\lambda(n)$ are $(1+1/h)$-dense in $[h,p_1\cdots p_J y]$. Since $p_J > Y^{(5/4)^{J-1}} \ge x^c$, this concludes the proof.
\end{proof}

\begin{proof}[Proof of Theorems \ref{dense} and \ref{denseh}]
Applying Theorem \ref{verydense}, there is a positive integer $k$ so that when $z$ is large enough, for more than half of the positive integers $d\le z$, $\phi(d)$ and $\lambda(d)$ are $(1+1/h)$-dense in $[h,z^{1/k}]$. Put $\eps = \frac{1}{5k^2}$, let $x$ be sufficiently large and $x^{\frac12} < d \le x^{\frac12+\eps}$, where $\phi(d)$ is $(1+1/h)$-dense in $[h, x^{\frac{1}{2k}}]$. Consider distinct primes $p_1$, $p_2, \ldots, p_k \in I := [x^{\frac{1}{2k}-2\eps}, x^{\frac{1}{2k}-\eps}]$ which do not divide $d$. Note that
\be\label{Dpi}
x^{1-2k\eps} \le d p_1 p_2 \cdots p_{k} \le x^{1-(k-1)\eps}.
\ee
Let $q$ be a prime not dividing $dp_1\cdots p_k$ and satisfying
\be\label{qrange1}
\frac12\frac{x}{dp_1 \cdots p_{k}} < q \le \frac{x}{d p_1 \cdots p_{k}},
\ee
so that by \eqref{Dpi} and the definition of $\eps$,
\be\label{qrange2}
x^{\frac1{6k}} \le q \le x^{\frac{2}{5k}}.
\ee

We claim that for all such numbers $n = d p_1 \cdots p_{k}q$ satisfying the additional hypothesis
\be\label{lambdabig}
\lambda(n) \ge x^{1-\eps},
\ee
$\phi(n)$ is $(1+1/h)$-dense in $[h,\phi(n)/(h + 1))$ and $\lambda(n)$ is $(1+1/h)$-dense in $[h,\lambda(n)/(h+1))$. Let $y=x^{\frac{1}{2k}}$. Observe that $\phi(n)=\phi(d) (p_1-1)\cdots (p_{k}-1)(q-1)$, $\phi(d)$ is $(1+1/h)$-dense in $[h,y]$, $p_i-1\le x^{\frac{1}{2k}-\eps} < (y/h)$ $(1\le i\le k)$ and $q\le (y/h)$. By Lemma \ref{densecondinterval} with $D=\phi(d)$, $m_i=p_i-1$ $(1\le i\le k)$ and $m_{k+1}=q-1$, $\phi(n)$ is $(1+1/h)$-dense in $[h,w]$, where $w=y(p_1-1)\cdots (p_k-1)(q-1)$. By \eqref{Dpi} and \eqref{qrange1}, 
$$
w \ge 2^{-k-1} y p_1\cdots p_k q 
\ge 2^{-k-2} \frac{x^{1+\frac{1}{2k}}}{d} \ge h \sqrt{x}.
$$
But $\phi(n)$ is also $(1+1/h)$-dense in $[\phi(n)/w,\phi(n)/(h + 1))$ since $d|m \iff (m/d)|m$, consequently $\phi(n)$ is $(1+1/h)$-dense in $[h, \phi(n)/(h + 1))$.

The argument for $\lambda(n)$ is similar, except that now
$$
\lambda(n)=\lambda(d) \frac{q-1}{f} \prod_{i=1}^k \frac{p_i-1}{f_i},
$$
where $f$ is some divisor of $q-1$ and $f_i$ is some divisor of $p_i-1$ $(1\le i\le k$). Here we use \eqref{lambdabig}, which implies that $f f_1\cdots f_k \le x^{\eps}$. By Lemma \ref{densecondinterval} with $D=\lambda(d)$, $m_i=(p_i-1)/f_i$ $(1\le i\le k)$ and $m_{k+1}=(q-1)/f$, we see that $\lambda(n)$ is $(1+1/h)$-dense in $[h,w]$, where 
$$
w=y \frac{q-1}{f} \prod_{i=1}^k \frac{p_i-1}{f_i} \ge
2^{-k-2} \frac{x^{1+\frac{1}{2k}-\eps}}{d} \ge h\sqrt{x}. 
$$
As with $\phi(n)$, we conclude that $\lambda(n)$ is $(1+1/h)$-dense in $[h, \lambda(n)/(h + 1))$.

Notice that for the above $n$, when $h = 1$, $\phi(n)$ is $2$-dense in $[1, \phi(n)/2)$. Since $\phi(n)$ is a divisor of itself, we conclude that $\phi(n)$ is $2$-dense in $[1, \phi(n))$ and hence $2$-dense. This conclusion also holds for $\lambda(n)$ by similar arguments.

Finally, we show that the number of such integers $n\le x$ is $\gg_h x$. First, \eqref{lambdabig} holds for almost all $n$ by Theorem 2 of \cite{EPS}. By the prime number theorem and \eqref{qrange2}, given $d,p_1,\ldots,p_{k}$, the number of possible primes $q$ is $\gg_k x/(dp_1\cdots p_{k} \log x)$. We also have
$$
\sum_{p_1,\ldots,p_k\in I} \frac{1}{p_1\cdots p_{k}} \gg_k 1,
$$
and $\sum 1/d \gg \log x$ by partial summation. Hence, there are $\gg_k x$
tuples $(d,p_1,\ldots,p_{k},q)$ with product $n\in (x/2,x]$ and with $\phi(n)$ and $\lambda(n)$ being $(1+1/h)$-dense respectively. Given such an integer $n$, $n$ has at most $6k$ prime factors $\ge x^{\frac{1}{6k}}$, hence the number of tuples $(d,p_1,\ldots,p_k,q)$ with product $n$ is bounded by a function of $k$. Thus the proof is complete.
\end{proof}

\begin{proof}[Proof of Theorem \ref{nondense}]
Suppose $0<c < 2\th-1$, and Let $\eps>0$ be so small that $2\th-1-6\eps>c$. Consider $n=pm\le x$, where $x^{1-2\eps} < p \le x^{1-\eps}$, and $P^+(p-1) > p^{\th-\eps}$. By the definition of $\th$, there are $\gg z/\log z$ such primes $\le z$, if $z$ is large enough. Then $\phi(n)$ and $\lambda(n)$ are each divisible by a prime $q$ with $q > x^{(1-2\eps)(\th-\eps)} > x^{\th-3\eps}$, and therefore neither function has divisors in $[x^{1-\th+3\eps},x^{\th-3\eps}]$. The number of such $n$ is, by partial summation,
$$
= \sum_{\substack{x^{1-2\eps}<p\le x^{1-\eps} \\ P^+(p-1)>p^{\th-\eps}}}
\left\lfloor \frac{x}{p} \right\rfloor \gg_{\eps} x,
$$
and the proof is complete.
\end{proof}

\begin{proof}[Proof of Corollary \ref{limitverydense}]
Let $\delta>0$. By Theorem \ref{verydense}, if $x$ is sufficiently large, then for at least $(1-\del)x$ integers $n\le x$, both $\phi(n)$ and $\lambda(n)$ are $(1+1/h)$-dense in $[h,x^{g(x)}]$. Since $\del$ is arbitrary, the corollary follows.
\end{proof}

\begin{proof}[Proof of Corollary \ref{bxy2y}]
(i) The elementary inequality $\sum_{n\le x} n/\phi(n) \ll x$ implies
that 
$$
|\{ n\le x : \phi(n) \le \eps n \}| \ll \eps x \qquad (0<\eps\le 1).
$$
Consequently, using Theorem \ref{dense}, if $c$ is small enough
then there are $\gg x$
of the integers $n\le x$ for which $\phi(n)$ is $2$-dense and
$\phi(n) \ge cx$. This proves (i) for $y\le cx$. For
a given constant $f\in [c,1/2]$, it is an elementary fact that
$fx < \phi(n) \le 2fx$ for $\gg_f x$ integers $n\le x$.
This completes the proof for the remaining $y$.

(ii) From the proof of Theorem \ref{nondense}, for a positive proportion
of integers $n$, $\phi(n)$ has no divisors in $[x^c,x^{1-c}]$.

(iii) This follows immediately from Corollary \ref{limitverydense}.
\end{proof}
\bibliographystyle{amsplain}
\bibliography{bxyz}

\end{document}